\documentclass[11pt]{amsart}
\usepackage{AFBoix2015}

\begin{document}

\title[On the infinitely generated locus of Frobenius algebras]{On the infinitely generated locus of Frobenius algebras of rings of prime characteristics$$\text{To our friend Ngo Viet Trung for his 70th birthday}$$}
\author[A.\,F.\,Boix]{Alberto F.\,Boix}
\address{IMUVA--Mathematics Research Institute, Universidad de Valladolid, Paseo de Belen, s/n, 47011, Valladolid, Spain.}
\email{alberto.fernandez.boix@uva.es}

\author[D.\,A.\,J.\,G\'omez--Ram\'irez]{Danny A.\,J.\,G\'omez--Ram\'irez}

\address{Parque Tech at the Instituci\'on Universitaria Pascual Bravo, Medell\'in, Colombia}
\email{daj.gomezramirez@gmail.com}

\author[S.\,Zarzuela]{Santiago Zarzuela}

\address{Departament de Matem\`atiques i Inform\`atica, Universitat de Barcelona, Gran Via de les Corts Catalanes 585, Barcelona 08007, SPAIN}
\email{szarzuela@ub.edu}

\keywords{Frobenius algebras, Stanley--Reisner rings.}

\subjclass[2020]{13A35, 13F55.}

\begin{abstract}
Let $R$ be a commutative Noetherian ring of prime characteristic $p$. The main goal of this paper is to study in some detail when
\[
\overline{W^R}:=\{\mathfrak{p}\in\operatorname{Spec} (R):\ \mathcal{F}^{E_{\mathfrak{p}}}\text{ is finitely generated as a ring over its degree zero piece}\}
\]
is an open set in the Zariski topology, where $\mathcal{F}^{E_{\mathfrak{p}}}$ denotes the Frobenius algebra attached to the injective hull of the residue field of $R_{\mathfrak{p}}.$ We show that this is true when $R$ is a Stanley--Reisner ring; moreover, in this case, we explicitly compute its closed complement, providing an algorithmic method for doing so. 
\end{abstract}

\maketitle

\section{Introduction}\label{Introduction}

The main goal of the present paper is to prove the following result: 

\begin{teo}

Let $S=\K[x_1, \dots, x_n]/I$ be a Stanley-Reisner ring, where $\K$ is a field of prime characteristic $p$ and $I$ a squarefree monomial ideal. Then, the set 
\[\overline{W^S} = \{\fp \in \operatorname{Spec} (S):\ \mathcal{F}^{E_{S_{\mathfrak{p}}}}\text{ is finitely generated as a ring over its degree zero piece}\},
\]where $\mathcal{F}^{E_{S_{\mathfrak{p}}}}$ is the Frobenius algebra of the injective hull of the residue field of the local ring $S_{\mathfrak{p}}$, is an open set of $\operatorname{Spec} (S)$. 

\end{teo}

Moreover, we determine the defining ideal of its closed complement, and provide an algorithmic method to compute it, both algebraic and combinatorial.    

\medskip 

We explain now the context and the meaning of the above statement. 

\medskip 

Let $R$ be a commutative Noetherian ring of prime characteristic $p,$ and let $M$ be an $R$--module. For each integer $e\geq 0,$ we denote by $\operatorname{End}_{p^e}(M)$ the set of $p^e$--linear maps of $M;$ that is, $\operatorname{End}_{p^e}(M)$ is made up by abelian group endomorphisms $\xymatrix@1{M\ar[r]^{\phi}& M}$ such that
\[
\phi (rm)=r^{p^e}\phi (m)\text{ for all }(r,m)\in R\times M.
\]
In this way, one can cook up the so--called Frobenius algebra of $M$
\[
\mathcal{F}^M:=\bigoplus_{e\geq 0}\operatorname{End}_{p^e}(M),
\]
where multiplication is given by composition of maps. This algebra, introduced in \cite{LyubeznikSmith2001} in the context of tight closure theory has attracted some attention in the last years.

One question raised in \cite[pp. 3156]{LyubeznikSmith2001} is whether this algebra is finitely generated over its degree zero piece. We want to briefly summarize here some of the answers that have appeared in the last years.

\begin{enumerate}[(i)]

\item If $M=R$ and $\xymatrix@1{R\ar[r]^-{F}& R}$ denotes the Frobenius endomorphism of $R,$ then $\mathcal{F}^R\cong R[\Theta;F],$ so the algebra is finitely generated \cite[Example 3.6]{LyubeznikSmith2001}. The same conclusion holds when $(R,\mathfrak{m})$ is a local complete $S_2$--ring, $M=H_{\mathfrak{m}}^{\dim (R)}(R)$ is the top local cohomology module of $R$ supported on $\mathfrak{m},$ and $F$ denotes the natural Frobenius action on this module \cite[Example 3.7]{LyubeznikSmith2001}.

\item If $M$ is an $R$--module of finite length, then $\mathcal{F}^M$ is finitely generated \cite[Theorem 2.11]{EnescuYao2014}.

\item Let $R=\K[\![x_1,\ldots ,x_n]\!]/I,$ where $\K$ is a field of prime characteristic $p,$ and $I$ is a squarefree monomial ideal. If $M=E_R$ is the injective hull of the residue field of $R,$ then either $\mathcal{F}^M$ is either infinitely generated or $\mathcal{F}^M\cong R[u\theta;F]$ for some $u\in R$ \cite[Theorem 3.5]{AlvarezBoixZarzuela2012}.

\item If $R$ is a normal, excellent, $\mathbb{Q}$--Gorenstein local ring of prime characteristic $p$ and order $m,$ and again $M=E_R$ is the injective hull of the residue field of $R,$ then $\mathcal{F}^M$ is finitely generated if and only if $\gcd(p,m)=1,$ see \cite[Proposition 4.1]{KatzmanSchwedeSinghZhang2013} and \cite[Theorem 4.5]{EnescuYao2014}.
\end{enumerate}

However, not too much is known about the following:

\begin{quo}\label{question about the openness of the locus}
Let $R$ be a commutative Noetherian ring of prime characteristic $p$. Is it true that
\[
\overline{W^R}:=\{\mathfrak{p}\in\operatorname{Spec} (R):\ \mathcal{F}^{E_{\mathfrak{p}}}\text{ is finitely generated as a ring over its degree zero piece}\},
\]
where $E_{\mathfrak{p}}$ is the injective hull of the residue field of $R_{\mathfrak{p}}$,
is an open set?
\end{quo}

In \cite{Gallegoetal2021}, it is shown that this question has a positive answer in some particular cases; however, to the best of our knowledge no systematic study of this problem seems to have been carried out so far.

The goal of this paper is to study Question \ref{question about the openness of the locus} in some detail; more precisely, let $S=R/I$ a commutative Noetherian ring which is a quotient of a regular Noetherian ring $R$ (that is, locally regular) of prime characteristic $p$. Under these assumptions, we show that $\overline{W^S}$ is closed under generalization; moreover, when $S$ is a Stanley--Reisner ring, we will show that $\overline{W^S}$ is really an open set in the Zariski topology, and provide an algorithmic and explicit description of the defining ideal of its closed complement.

Our approach to the problem is by means of the study of a similar question over another non-commutative graded algebra, that we denote by $\mathcal{A}^S$ (see Definition \ref{non principally generated locus introduced}). This algebra is much easier to handle and behaves well under flat morphisms. Moreover, thanks to a well known result of R. Fedder \cite{Fedder1983}, if $R$ is a complete regular local ring, then the algebra $\mathcal{A}^S$ is isomorphic to $\mathcal{F}^E$, where $E$ is the injective envelope of the residue field of $S$, see Section \ref{background} for the details.

It turns out that, in the Stanley-Reisner case, the graded pieces of the algebra $\mathcal{A}^S$ are homogeneous for the natural multigrading of the polynomial ring $R = \K[x_1, \dots , x_n]$. This allows to control their local vanishing, and in the end the finite generation of $\mathcal{A}^{S_{\fp}}$, just by looking at the localization by homogeneous prime ideals of $R$, which are a finite set. Analyzing carefully this process, we are able to prove that the complement of the set $\overline{W^S}$ is the closed set defined by a square free monomial ideal, see Theorem \ref{the infinitely generated locus is closed}. 

In order to find explicitly the above defining ideal, first we need to pass from the local case to the Stanley-Reisner case. This, we do by relating through a faithfully flat morphism the localization of $S$ at an homogeneous prime ideal to an adequate monomial localization (see Definition \ref{introducing monomial localization}), which provides the needed frame to apply the criteria from \cite{AlvarezBoixZarzuela2012} and \cite{AlvarezYanagawa2014} to determine whether the corresponding Frobenius algebra of the injective hull of the residue field is finitely generated or not in the Stanley-Reisner case. Then, we introduce a very simple algebraic algorithm that computes explicitly our defining ideal, see Theorem \ref{main result: calculation of the infinitely generated locus}. From the combinatorial point of view, this process corresponds to consider all the possible links of the simplicial complex determined by $S$, and to look for the existence of the so named free faces inside them. The details are developed in Section \ref{algorithm}.

Finally, in Section \ref{Examples} we provide several examples computed by means of the implentation of our algorithm in Macaulay2.

\section{The non--finitely generated locus: definition and general properties}\label{background}

\begin{df}\label{non principally generated locus introduced}
Let $\K$ be a field of prime characteristic $p>0,$ let $R$ be a commutative Noetherian ring containing $\K,$ $S=R/I,$ where $I\subset R$ is an ideal, and set
\[
\mathcal{A}^S:=\bigoplus_{e\geq 0}\left(\frac{(I^{[p^e]}:_R I)}{I^{[p^e]}}\right)=\bigoplus_{e\geq 0}\mathcal{A}_e.
\]
This is a graded abelian group that we can endow with a non--commutative ring structure by setting its multiplication as
\[
a\cdot b=ab^{p^e},\ a\in\mathcal{A}_e,\ b\in\mathcal{A}_{e'}.
\]
Finally, set
\[
\mathcal{A}^{S_{\mathfrak{p}}}:=\bigoplus_{e\geq 0}\left(\frac{(I_{\mathfrak{p}}^{[p^e]}:_{R_{\mathfrak{p}}} I_{\mathfrak{p}})}{I_{\mathfrak{p}}^{[p^e]}}\right),
\]
where again this is regarded as a non--commutative ring with multiplication defined as before. Then, we refer to
\[
W:=\{\mathfrak{p}\in\operatorname{Spec}(R):\ \mathfrak{p}\supseteq I,\ \mathcal{A}^{S_{\mathfrak{p}}}\text{ is a not finitely generated ring over }S_{\mathfrak{p}}\},
\]
as the \textbf{non--finitely generated locus} of the algebra $\mathcal{A}^S.$ Moreover, we denote by $\overline{W}$ the complement of $W$ inside $V(I).$
\end{df}
The reader might ask why to care about the non--commutative algebra defined above; actually, this ring is very close to the Frobenius algebra attached to the injective hull of the residue field of a complete local ring, as the below discussion explains.

\begin{disc}\label{preliminarias of completion to correctly talk about localization}
Let $(R,\mathfrak{m})$ be a commutative Noetherian regular local ring of prime characteristic $p,$ let $I\subset R$ be an ideal, and set $S:=R/I.$ Moreover, let $T:=\widehat{S}$ be the completion of $S$ with respect to the $\mathfrak{m}$--adic topology. Now, consider
\[
\mathcal{G}^S=\bigoplus_{e\geq 0}\left(\frac{(I^{[p^e]}:_R I)}{I^{[p^e]}}\right)F^e=\bigoplus_{e\geq 0}\mathcal{G}_e F^e,
\]
where $F^e$ denotes the $e$--th iteration of the Frobenius map on the injective hull of the residue field of $R.$ It is clear that $\mathcal{G}^S$ is an $\N_0$--graded ring, not necessarily commutative with $\mathcal{G}_0=S$ which is degree--wise finitely generated and a left $S$--skew algebra; now, we also consider
\[
\mathcal{G}^T:=\bigoplus_{e\geq 0}\mathcal{B}_e F_T^e,\text{ where }\mathcal{B}_e:=T\otimes_S \mathcal{G}_e,
\]
and $F_T^e$ denotes the $e$--th iteration of the Frobenius map on $T.$
However, notice that, by Fedder's Lemma \cite[pp. 465]{Fedder1983}, $\mathcal{G}^T$ is isomorphic to $\mathcal{F}^{E_S},$ the Frobenius algebra attached to the injective hull of the residue field of $S,$ which in turn is also isomorphic to $\mathcal{F}^{E_T}$ because of \cite[Proposition 3.3]{LyubeznikSmith2001}. On the other hand, we also have a graded ring isomorphism
\[
\mathcal{A}^S\cong\mathcal{G}^S.
\]
These facts allow us to regard $\mathcal{A}^S$ as a subring of $\mathcal{F}^{E_T}.$
\end{disc}


Before going on, we need to review the following notions borrowed from \cite{EnescuYao2014}.

\begin{disc}\label{reminders about complexity sequence}
Let $\cA=\bigoplus_{e\geq 0}\cA_e$ be an $\N_0$--graded ring, not necessarily commutative.

\begin{enumerate}[(i)]

\item Let $G_e:=G_e (\cA)$ be the subring of $\cA$ generated by the homogeneous elements of degree
less than or equal than $e;$ we agree that $G_{-1}=\cA_0.$ Notice also that $G_0=\cA_0.$

\item Let $k_e:=k_e (\cA)$ be the minimal number of homogeneous generators of $G_e$ as a subring of
$\cA$ over $\cA_0;$ we agree that $k_{-1}=0.$ One says that $\cA$ is \textbf{degree--wise finitely
generated} if $k_e<\infty$ for any $e.$


\item Again, assume that $\cA$ is a degree--wise finitely generated ring and set, for each
integer $e\geq 0,$ $c_e (\cA):=k_e-k_{e-1};$ as pointed out in \cite[Remark 2.6]{EnescuYao2014}, $\cA$
is finitely generated as a ring over $\cA_0$ if and only if the sequence
$\{c_e(\cA)\}_{e\geq 0}$ is eventually zero.

\item Let $R$ be a commutative ring, and let $\cA$ be a degree--wise finitely generated ring such that
$R=\cA_0;$ moreover, assume that $\cA$ is a left \textbf{$R$--skew algebra} (i.e.
$aR\subseteq Ra$ for all homogeneous elements $a\in\cA$). Then, by \cite[Corollary 2.10]{EnescuYao2014}
$c_e (\cA)$ equals the minimum number of generators of $\cA_e/(G_{e-1})_e$ as a left $R$--module
for any $e.$

\end{enumerate}

\end{disc}
As observed along Discussion \ref{reminders about complexity sequence}, we know that, when $\cA$ is a non--commutative, degree--wise finitely generated ring, $\cA$ is finitely generated as a ring over its degree zero piece if and only if the sequence
$\{c_e(\cA)\}_{e\geq 0}$ is eventually zero. This fact motivates us to introduce the following:

\begin{df}\label{k-generation}
Given $\cA$ a non--commutative, degree--wise finitely generated ring, and given an integer $k\geq 1,$ we say that $\cA$ is \textbf{k--generated} provided $c_e(\cA)=0$ for all $e>k;$ when $k=1,$ we say that $\cA$ is \textbf{principally generated}.
\end{df}

Now, we are in position to introduce the following:

\begin{df}\label{stratifying the infinitely generated locus}
Under the assumptions and notations of Definition \ref{non principally generated locus introduced}, we can define
\[
\overline{W}_k:=\{\mathfrak{p}\in\operatorname{Spec}(R):\ \mathfrak{p}\supseteq I,\ \cA^{S_{\mathfrak{p}}}\text{ is $k$--generated}\}.
\]
The reader will easily note that the $\overline{W}_k$'s provide a stratification of the set
\[
\overline{W}=\{\mathfrak{p}\in\operatorname{Spec}(R):\ \mathfrak{p}\supseteq I,\ \cA^{S_{\mathfrak{p}}}\text{ is finitely generated over }S_{\mathfrak{p}}\};
\]
in other words,
\[
\overline{W}=\bigcup_{k\geq 1}\overline{W}_k.
\]
\end{df}

Our next goal is to compute the sequence $\{c_e(\cA)\}_{e\geq 0}$ for the ring $\mathcal{A}^S,$ and to characterize when it is eventually zero.

\begin{prop}\label{characterization of eventually zero sequence}
Let $R$ be a commutative Noetherian ring of prime characteristic $p,$ let $I\subset R$ be an ideal, and set $S:=R/I.$ Moreover, for each $e\geq 1$ write $K_e:=(I^{[p^e]}:_R I),$ and set
\[
L_e:=\sum_{\substack{1\leq a_1,\ldots,\ a_s\leq e-1\\ a_1+\ldots+a_s=e}}K_{a_1}K_{a_2}^{[p^{a_1}]}K_{a_3}^{[p^{a_1+a_2}]}\cdots K_{a_s}^{[p^{a_1+\ldots+a_{s-1}}]}.
\]
Then, the following assertions hold.

\begin{enumerate}[(i)]

\item For any $e\geq 1,$ one has that $c_e(\mathcal{A}^S)$ equals the minimum number of generators as a left $S$--module of
\[
\frac{(I^{[p^e]}:_R I)}{L_e}.
\]

\item The sequence $\{c_e(\mathcal{A}^S)\}_{e\geq 0}$ is eventually zero if and only if, for all $e\gg 0,$ one has that
\[
(I^{[p^e]}:_R I)=L_e.
\]

\item Given an integer $k\geq 1,$ the ring $\cA^S$ is $k$--generated if and only if, for all $e>k,$ one has that
\[
(I^{[p^e]}:_R I)=L_e.
\]

\end{enumerate}

\end{prop}

\begin{proof}
Since parts (ii) and (iii) follow immediately from part (i), we only plan to prove part (i); indeed, as viewed in Discussion \ref{reminders about complexity sequence}, $c_e(\mathcal{A}^S)$ equals the minimum number of generators as a left $S$--module of $\cA_e/(G_{e-1}\cap\cA_e),$ where $G_{e-1}$ is as defined in Discussion \ref{reminders about complexity sequence}. However, thanks to \cite[Proposition 2.1]{Katzman2010} we know that
\[
G_{e-1}\cap\cA_e=L_e.
\]
In this way, the result follows directly from this fact.
\end{proof}

\begin{rk}
Notice that \cite[Proposition 2.1]{Katzman2010} was only proved for a formal power series ring in three indeterminates; however, the reader can easily verify that its proof also holds in our setting.
\end{rk}
Now, our aim is to show how behaves the algebra $\cA^S$ under certain base change; this behaviour will play a key role in several places of this paper later.

\begin{prop}\label{quasi Frobenius algebra and base change}
Let $\xymatrix@1{R\ar[r]^-{\phi}& R'}$ be a faithfully flat ring homomorphism between commutative Noetherian regular rings of prime characteristic $p,$ let $I\subset R$ be an ideal, set $I':=\phi (I)R'$ and $S':=R'/I'.$ Now, set
\[
\cA^S:=\bigoplus_{e\geq 0} \left(\frac{(I^{[p^e]}:_R I)}{I^{[p^e]}}\right),\quad \cA^{S'}=\bigoplus_{e\geq 0} \left(\frac{(I'^{[p^e]}:_{R'} I')}{I'^{[p^e]}}\right).
\]

Finally, we assume, for any integer $e\geq 0,$ that
\[
\phi ((I^{[p^e]}:_R I))R'=(I'^{[p^e]}:_{R'} I').
\]
Then, $\cA^S$ is $k$--generated for some integer $k\geq 1$ if and only if $\cA^{S'}$ is $k$--generated; in particular, this implies that $\cA$ is finitely generated as a ring over $S,$ if and only if $\cA^{S'}$ is finitely generated as a ring over $S'.$





\end{prop}

\begin{proof}
We know that $\cA^S$ is $k$--generated as a ring over $S$ if and only if $c_e (\cA^S)=0$ for all $e>k$, which is equivalent to say, thanks to Proposition \ref{characterization of eventually zero sequence}, that for all $e> k,$ one has $(I^{[p^e]}:_R I)=L_e,$ where
\[
L_e=\sum_{\substack{1\leq a_1,\ldots, a_s\leq e-1\\ a_1+\ldots+a_s=e}}(I^{[p^{a_1}]}:_R I)\cdot (I^{[p^{a_2}]}:_R I)^{[p^{a_1}]}\cdots (I^{[p^{a_s}]}:_R I)^{[p^{a_1+\ldots+a_{s-1}}]}.
\]
In this way, if for all $e>k,$ one has that $(I^{[p^e]}:_R I)=L_e,$ then this is equivalent to say, thanks to our assumptions, that for all $e>k$ $(I'^{[p^e]}:_{R'} I')=L'_e,$ where
\[
L'_e=\sum_{\substack{1\leq a_1,\ldots, a_s\leq e-1\\ a_1+\ldots+a_s=e}}(I'^{[p^{a_1}]}:_{R'} I')\cdot (I'^{[p^{a_2}]}:_{R'} I')^{[p^{a_1}]}\cdots (I'^{[p^{a_s}]}:_{R'} I')^{[p^{a_1+\ldots+a_{s-1}}]}.
\]
The proof is therefore complete.
\end{proof}
By applying Proposition \ref{quasi Frobenius algebra and base change} to the completion map, we immediately obtain the following:

\begin{cor}\label{we can restrict our attention to the algebra A}
Let $(R,\mathfrak{m})$ be a commutative Noetherian regular local ring of prime characteristic $p,$ let $I\subset R$ be an ideal, and set $S:=R/I.$ Given an integer $k\geq 1,$ one has that $\cA^S$ is $k$--generated if and only if $\mathcal{F}^{E_S}$ is $k$--generated; in particular, $\cA^S$ is finitely generated as a ring over its degree zero piece if and only if so is $\mathcal{F}^{E_S}.$
\end{cor}
Our next goal is to show that the set $W$ introduced in Definition \ref{non principally generated locus introduced} is closed under specialization. This is exactly the content of the next:

\begin{prop}\label{the principal locus is stable under generalization}
Let $k\geq 1$ be an integer, and let
\[
\overline{W}_k=\{\mathfrak{p}\in\operatorname{Spec}(R):\ \mathfrak{p}\supseteq I,\ \mathcal{A}^{S_{\mathfrak{p}}}\text{ is $k$--generated}\}.
\]
Then, $\overline{W}_k$ is closed under generalization. Equivalently, given prime ideals $\mathfrak{q}\subset\mathfrak{p}$ with $\mathfrak{p}\in\overline{W}_k,$ one has that $\mathfrak{q}\in\overline{W}_k.$

Therefore, $W_k$ is closed under specialization and therefore by \cite[\href{https://stacks.math.columbia.edu/tag/0EES}{Tag 0EES}]{stacks-project}, $W_k$ can be expressed as a directed union of closed subsets of $V(I).$

In particular, the set
\[
\overline{W}=\{\mathfrak{p}\in\operatorname{Spec}(R):\ \mathfrak{p}\supseteq I,\ \mathcal{A}^{S_{\mathfrak{p}}}\text{ is a finitely generated ring over }S_{\mathfrak{p}}\}.
\]
is closed under generalization.
\end{prop}

\begin{proof}
Let $\mathfrak{q}\subset\mathfrak{p}$ be prime ideals with $\mathfrak{p}\in\overline{W},$ we assume that $\mathcal{A}^S$ is not $k$-- generated, otherwise the statement is obvious (indeed, in that case, $W_k=\emptyset.$) Let $\xymatrix@1{R\ar[r]^-{l_{\mathfrak{p}}}& R_{\mathfrak{p}}\ar[r]^-{l_{\mathfrak{p},\mathfrak{q}}}& R_{\mathfrak{q}}}$ be the natural localization maps, notice that $l_{\mathfrak{q}}=l_{\mathfrak{p},\mathfrak{q}}\circ l_{\mathfrak{p}}.$ Since $\mathfrak{p}\in\overline{W},$ one has that $c_e(\cA^{S_{\mathfrak{p}}})=0$ for all $e>k,$ which means, due to Proposition \ref{characterization of eventually zero sequence}, that for all $e>k,$
\[
(I^{[p^e]}:_R I)R_{\mathfrak{p}}=L_e R_{\mathfrak{p}}.
\]
Applying to this equality the map $l_{\mathfrak{p},\mathfrak{q}}$ and using the equality $l_{\mathfrak{q}}=l_{\mathfrak{p},\mathfrak{q}}\circ l_{\mathfrak{p}}$ one finally obtains that for all $e>k,$
\[
(I^{[p^e]}:_R I)R_{\mathfrak{q}}=L_e R_{\mathfrak{q}}.
\]
In this way, this is equivalent to say, using once again Proposition \ref{characterization of eventually zero sequence}, that $c_e(
\cA^{S_{\mathfrak{q}}})=0$ for all $e>k$, and therefore $\mathfrak{q}\in\overline{W}_k,$ as claimed.
\end{proof}

As immediate consequence of Proposition \ref{the principal locus is stable under generalization} one gets the below:

\begin{cor}\label{the locus of Frobenius algebras is stable under generalization}
Let $R$ be a commutative Noetherian regular ring of prime characteristic $p,$ let $I\subset R$ be an ideal, and set $S:=R/I.$ Then, for any integer $k\geq 1,$
\[
\{\mathfrak{p}\in\operatorname{Spec}(R):\ \mathfrak{p}\supseteq I,\ \mathcal{F}^{E_{S_{\mathfrak{p}}}}\text{ is $k$--generated}\}
\]
is closed under generalization; in particular,
\[
\{\mathfrak{p}\in\operatorname{Spec}(R):\ \mathfrak{p}\supseteq I,\ \mathcal{F}^{E_{S_{\mathfrak{p}}}}\text{ is a finitely generated ring over }\widehat{S_{\mathfrak{p}}}\}
\]
is closed under generalization.
\end{cor}

\section{The case of a Stanley--Reisner ring}\label{section of preliminaries}
In what follows in this section, let $\K$ be a field, and let $R=\K [x_1,\ldots ,x_n]$ be a polynomial ring over $\K.$ We regard $R$ as an $\mathbb{N}^n$--graded ring, where $\deg (x_i)$ is the $i$--th canonical basis vector in $\mathbb{N}^n.$ With this grading, the graded $R$--submodules of $R$ (aka the graded ideals of $R$) are exactly the monomial ideals, and the graded prime ideals of $R$ are the ideals generated by a subset of the variables; in this way, we denote by $\hbox{}^*\operatorname{Spec} (R)$ the set of graded prime ideals of $R.$

In this setting, given $\mathfrak{p}$ a (not necessarily graded) prime ideal of $R$ containing a monomial ideal $J\subset R,$ it is known that $\mathfrak{p}\supseteq\mathfrak{p}^*\supseteq J,$ where $\mathfrak{p}^*$ is the graded prime ideal of $R$ generated by the homogeneous elements of $\mathfrak{p}.$

The first technical result we need to establish is that, in order to calculate the non--finitely generated locus of a Stanley--Reisner ring, it is enough to restrict to face ideals; this holds because of the following:

\begin{lm}\label{check only with graded prime}
Let $I\subset R$ be a squarefree monomial ideal, set $S:=R/I,$ and assume that there is a prime $\mathfrak{p}\in\operatorname{Spec} (R)$ such that $\mathcal{A}^{S_{\mathfrak{p}^*}}$ is $k$--generated for some integer $k\geq 1.$ Then, $\mathcal{A}^{S_{\mathfrak{p}}}$ is also $k$--generated.

In particular, if there is a prime $\mathfrak{p}\in\operatorname{Spec} (R)$ such that $\mathcal{A}^{S_{\mathfrak{p}^*}}$ is finitely generated as a ring over $S_{\mathfrak{p}^*},$ then we have that $\mathcal{A}^{S_{\mathfrak{p}}}$ is also finitely generated as a ring over $S_{\mathfrak{p}}.$
\end{lm}

\begin{proof}
Assume that $\mathcal{A}^{S_{\mathfrak{p}^*}}$ is $k$--generated as a ring over $S_{\mathfrak{p}^*};$ this means that $c_e(\mathcal{A}^{S_{\mathfrak{p}^*}})=0$ for all $e>k$, hence for all $e>k$ one has that
\[
\left(\frac{(I^{[p^e]}:_R I)}{L_e}\right)_{\mathfrak{p}^*}=0.
\]
This means, setting $M_e:=\frac{(I^{[p^e]}:_R I)}{L_e},$ that $\mathfrak{p}^*\notin\operatorname{Supp} (M_e),$ and this implies, by \cite[13.1.6 (i)]{BroSha}, that $\mathfrak{p}\notin\operatorname{Supp} (M_e)$ for all $e>k.$ But, this is equivalent to say that $c_e(\mathcal{A}^{S_{\mathfrak{p}}})=0$ for all $e>k$, which implies that $\mathcal{A}^{S_{\mathfrak{p}}}$ is $k$--generated; the proof is therefore complete.
\end{proof}
Now, we are in position to prove one of the main results of this paper; namely, the following:

\begin{teo}\label{the infinitely generated locus is closed}
Let $\K$ be a field of prime characteristic $p,$ let $I\subseteq\K [x_1,\ldots, x_n]=R$ be a squarefree monomial ideal, and set $S:=R/I.$ Then, the set
\[
W=\{\mathfrak{p}\in\operatorname{Spec}(R):\ \mathfrak{p}\supseteq I,\ \mathcal{A}^{S_{\mathfrak{p}}}\text{ is not finitely generated as a ring over }S_{\mathfrak{p}}\},
\]
is a closed set in the Zariski topology.
\end{teo}

\begin{proof}
If $W=\emptyset$ then we are done, so we can assume $W\neq\emptyset;$ now, let $\mathfrak{p}\in W.$ Thanks to Lemma \ref{check only with graded prime}, we have that $\mathfrak{p}^*\in W.$ This shows that the minimal members of the set $W$ are face ideals, which turn out to be a finite set. Thus, let $\mathfrak{p}_1^*,\ldots, \mathfrak{p}_t^*$ be the face ideals that belong to $W,$ and set
\[
J:=\bigcap_{a=1}^t \mathfrak{p}_a^*.
\]
Our above argument shows that $W\subseteq V(J);$ conversely, let $\mathfrak{q}\supset J,$ in particular $\mathfrak{q}\supset\mathfrak{p}_a$ for some $a=1,\ldots ,t.$ In this way, since $\mathfrak{p}_a\in W$ and $W$ is closed under specialization by Proposition \ref{the principal locus is stable under generalization}, we have that $\mathfrak{q}\in W.$

Summing up, we have finally checked that $W=V(J),$ hence a Zariski closed set.
\end{proof}
Since we know that, in the Stanley--Reisner case, the finite generation of the Frobenius algebra is equivalent to its principal generation (equivalently, to its $1$--generation using the terminology that we employ in this paper), Theorem \ref{the infinitely generated locus is closed} implies the following:

\begin{teo}\label{the principally generated locus is closed}
Let $\K$ be a field of prime characteristic $p,$ let $I\subseteq\K [x_1,\ldots, x_n]=R$ be a squarefree monomial ideal, and set $S:=R/I.$ Then, we have
\[
\overline{W}=\overline{W}_1=\{\mathfrak{p}\in\operatorname{Spec}(R):\ \mathfrak{p}\supseteq I,\ \mathcal{A}^{S_{\mathfrak{p}}}\text{ is $1$--generated}\}.
\]
\end{teo}

\begin{rk}
Theorem \ref{the infinitely generated locus is closed} is equivalent to say that
\[
\overline{W}=\{\mathfrak{p}\in\operatorname{Spec}(R):\ \mathfrak{p}\supseteq I,\ \mathcal{A}^{S_{\mathfrak{p}}}\text{ is finitely generated as a ring over }S_{\mathfrak{p}}\}
\]
is open in the Zariski topology. The reader will easily note that our proof does not involve the use of the so-called (the terminology is borrowed from \cite{Kimura2022}) topological Nagata criterion \cite[Theorem 24.2]{Mat86}.
\end{rk}

\begin{rk} Similarly as we have already done in Corollary \ref{the locus of Frobenius algebras is stable under generalization}, both Theorem \ref{the infinitely generated locus is closed} and Theorem \ref{the principally generated locus is closed} can also be equally formulated for the corresponding Frobenius algebra.
\end{rk}

\section{An algorithmic description of the non--finitely generated locus}\label{algorithm}

Our goal now is to give, in the case of a Stanley--Reisner ring, an explicit and algorithmic description of the non--finitely generated locus of its corresponding Frobenius algebra. We continue with the same notations as in Section \ref{section of preliminaries}.

Hereafter in this section, $[n]$ will denote the subset $\{1,\ldots, n\}$ of $n\geq 1$ elements, let $\Delta\subset [n]$ be a simplicial complex, and let $I=I_{\Delta}\subseteq R$ be the squarefree monomial ideal attached to $\Delta$ via the Stanley correspondence \cite[1.6 and 1.7]{MillerSturmfels2005}; moreover, given $F\subset\Delta$ a face, we denote by $\mathbf{p}_F$ the prime ideal of $R$ generated by the variables whose indices are not in $F.$ In other words,
\[
\mathfrak{p}_F :=(x_i:\ i\notin F).
\]
On the other hand, given a monomial $m=x_1^{a_1}\cdots x_n^{a_n},$ we denote by $\operatorname{supp} (m)$ its support; that is,
\[
\operatorname{supp}(m)=\{i\in [n]:\ a_i\neq 0\}.
\]
Finally, set
\[
\mathbf{x}_F:=\prod_{i\in F} x_i.
\]

\begin{lm}\label{a colon ideal in the middle}
If $F$ is a face of $\Delta,$ then $I\subseteq (I:_R \mathbf{x}_F)\subseteq\mathfrak{p}_F.$
\end{lm}

\begin{proof}
Since the inclusion $I\subseteq (I:_R \mathbf{x}_F)$ holds by definition of colon ideals, it is enough to check that $(I:_R \mathbf{x}_F)\subseteq\mathfrak{p}_F.$

Indeed, let $m_1,\ldots, m_t$ be a system of squarefree monomial generators of $I,$ it is known \cite[1.2.2]{HerzogHibi2011} that $(I:_R \mathbf{x}_F)=(g_1,\ldots, g_t),$ where
\[
g_j=\frac{m_j}{\gcd (m_j,\mathbf{x}_F)},\ 1\leq j\leq t,
\]
Therefore, it is enough to check that, given $1\leq j\leq t,$ $g_j$ is divisible by a variable $x_i$ with $i\notin F.$

So, fix $1\leq j\leq t,$ and set $S_j:=\operatorname{supp} (m_j);$ notice that, if $S_j\subset F,$ then $m_j$ would divide $\mathbf{x}_F$ and therefore $\mathbf{x}_F\in I,$ a contradiction because $F$ is a face of $\Delta.$ Hence there is $i\in S_j$ such that $i\notin F,$ thus $x_i$ divides, not only $m_j,$ but also $g_j,$ as claimed.
\end{proof}
The reader might ask why one needs to consider the above colon ideal. Next statement gives a combinatorial reason; recall that, given a simplicial complex $\Delta$ and a face $F$ of it, the link of $F$ inside $\Delta$ is defined as
$\link(F):=\{G\subset\Delta:\ F\cap G=\emptyset,\ F\cup G\subset\Delta\}.$

\begin{lm}\label{saturation and links}
Let $F\subset\Delta$ be a face. Then, $(I:_R \mathbf{x}_F)=I_{\link (F)}.$

\end{lm}

\begin{proof}
First of all, let $m_1,\ldots, m_t$ be a system of squarefree monomial generators of $I;$ since \cite[1.2.2]{HerzogHibi2011} $(I:_R \mathbf{x}_F)=(g_1,\ldots, g_t),$ where
\begin{equation}\label{expression of generators}
g_j=\frac{m_j}{\gcd (m_j,\mathbf{x}_F)},\ 1\leq j\leq t,
\end{equation}
one has that $(I:_R \mathbf{x}_F)$ is also a squarefree monomial ideal, so $(I:_R \mathbf{x}_F)=I_{\Delta '}$ for some simplicial complex $\Delta '.$ Moreover, as $I\subseteq (I:_R \mathbf{x}_F),$ one has, again as consequence of the Stanley--Reisner Correspondence, that $\Delta '\subset\Delta.$ Finally, given $G\subset [n]$ such that $\mathbf{x}_G$ is a squarefree minimal monomial generator of $(I:_R \mathbf{x}_F),$ one has that $G\notin\Delta '$ if and only if $\mathbf{x}_G\in (I:_R \mathbf{x}_F),$ which is equivalent to say that $\mathbf{x}_F\cdot\mathbf{x}_G\in I,$ which is equivalent to say (notice that $F\cap G=\emptyset$ because of \eqref{expression of generators}) that $\mathbf{x}_{F\cup G}\in I,$ which is equivalent to say that $F\cup G\notin\Delta.$
\end{proof}

Now, before establishing the main result of this section, we want to review for the convenience of the reader the notion of monomial localization as introduced in \cite[pp. 293]{HerzogRaufVladoiu2013}.

\begin{df}\label{introducing monomial localization}
Let $J\subseteq R$ be a (not necessarily squarefree) monomial ideal. We define the \textbf{monomial localization of $J$ at the prime ideal $\mathfrak{p}_F$} as the monomial ideal $J(\mathfrak{p}_F)\subset\K[x_i:\ i\notin F]$ obtained from $J$ by setting $x_j=1$ for all variables $j\in F.$ In other words, $J(\mathfrak{p}_F)$ is the extension of $J$ with respect to the $\K$--algebra map
\begin{align*}
\varphi_{\mathfrak{p}_F}:\ & R\rightarrow\K[x_i:\ i\notin F],\\
& x_j\longmapsto\begin{cases} x_j,\text{ if }j\notin F,\\ 1,\text{ if }j\in F. \end{cases}
\end{align*}
\end{df}

We also need to establish the below technical statement that will play a key role along the proof of the main result of this section. In the following result, given $J$ a monomial ideal with minimal monomial generating set $\{m_1,\ldots,m_r\},$ we denote by $\lcm_J$ the following monomial ideal:
\[
\lcm_J:=(\lcm(m_i,m_j):\ 1\leq i<j\leq r).
\]
Finally, given $I=(f_1,\ldots,f_s)$ an ideal inside a commutative Noetherian ring, we set $I^{[2]}:=(f_1^2,\ldots,f_s^2).$ 

\begin{lm}\label{links and localized Frobenius algebras}
Let $I\subset R$ be a squarefree monomial ideal, and let $\mathfrak{p}=\mathfrak{p}_F\in\hbox{}^*\operatorname{Spec} (R)$ be a face ideal. Then, the following statements are equivalent.

\begin{enumerate}[(i)]

\item $\mathcal{A}^{S_{\mathfrak{p}}}$ is a finitely generated ring over $S_{\mathfrak{p}}.$

\item $\mathcal{A}^{\widetilde{R}}$ is a finitely generated ring over $\widetilde{R},$ where $\widetilde{R}:=\K[x_i:\ i\notin F]_{(x_i:\ i\notin F)}.$

\item $\mathcal{A}^{R'}$ is a finitely generated ring over $R',$ where $R':=\K [x_i:\ i\notin F].$

\item We have $(I'^{[2]}:_{R'} I')=I'^{[2]}+(\lcm_{I'}),$ where $I':=I(\mathfrak{p}_F).$

\item We have $(K^{[2]}:_R K)= K^{[2]}+(\lcm_K),$ where $K:=(I:_R \mathbf{x}_F).$

\end{enumerate}

\end{lm}

\begin{proof}
We consider the following maps, where starting from the left, the first one is just $\varphi:=\varphi_{\mathbf{p}_F}$ as defined in Definition \ref{introducing monomial localization}, the second one is the obvious inclusion, and the last one is localization at $\mathfrak{p}.$
\[
\xymatrix{R\ar[r]^-{\varphi}& R'\ar@{^(->}[r]& R\ar[r]& R_{\mathfrak{p}}.}
\]
Since the inclusion $\xymatrix@1{R'\,\ar@{^(->}[r]& R}$ is given by adjoining variables and $\xymatrix@1{R\ar[r]& R_{\mathfrak{p}}}$ is a localization, the composition $\xymatrix@1{R'\ar[r]& R_{\mathfrak{p}}}$ is flat; moreover, the ideal $(x_i:\ i\notin F)\subset R'$ maps under this composition to $\mathfrak{p}R_{\mathfrak{p}}.$ Combining these two facts one gets that the induced inclusion
\[
\xymatrix{\widetilde{R}=\K[x_i:\ i\notin F]_{(x_i:\ i\notin F)}\ar[r]& R_{\mathfrak{p}}}
\]
is faithfully flat; this shows, combined with Proposition \ref{quasi Frobenius algebra and base change} that parts (i) and (ii) are equivalent.

Moreover, since $(x_i:\ i\notin F)\subset R'$ is the unique homogeneous maximal ideal of the graded ring $R',$ one has that the localization $\xymatrix@1{R'\ar[r]& \widetilde{R}}$, restricted to graded $R'$--modules, is also faithfully flat, hence parts (ii) and (iii) are also equivalent again using Proposition \ref{quasi Frobenius algebra and base change}.

Now, notice that $R'\otimes_R (R/I)\cong R'/I':=S'$ is also a Stanley--Reisner ring; therefore, by \cite[Theorem 3.5]{AlvarezBoixZarzuela2012} and \cite[Remark 2 and Lemma 3]{AlvarezYanagawa2014}, the algebra $\mathcal{A}^{S'}$ is finitely generated as ring over $S'$ if and only if
\[
(I'^{[2]}: I')= I'^{[2]}+(\lcm_{I'}),
\]
which proves the equivalence between parts (iii) and (iv).

Finally, set $K:=(I:_R \mathbf{x}_F).$ Along Lemma \ref{saturation and links} we saw that $K$ admits a minimal monomial generating set $G$ such that, if a monomial $m$ belongs to $G,$ then the support of $m$ is disjoint from $F.$ On the other hand, notice, by definition of $\varphi,$ that $\varphi(K)R'=I';$ actually, both $K$ and $I'$ admit $G$ as minimal generating set. This implies that
\[
(I'^{[2]}:_{R'} I')=I'^{[2]}+(\lcm_{I'}),
\]
if and only if
\[
(K^{[2]}:_R K)= K^{[2]}+(\lcm_K),
\]
just what we finally wanted to show.
\end{proof}

Finally, we are in position to establish the second main result of this paper, which is the following:

\begin{teo}\label{main result: calculation of the infinitely generated locus}
Let $\K$ be a field of prime characteristic $p,$ let $R=\K [x_1,\ldots, x_n],$ $S=R/I,$ and $I$ is a squarefree monomial ideal. Then, the set
\[
W=\{\mathfrak{p}\in\operatorname{Spec}(R):\ \mathfrak{p}\supseteq I,\ \mathcal{A}^{S_{\mathfrak{p}}}\text{ is not finitely generated as a ring over }S_{\mathfrak{p}}\}
\]
\[
=\{\mathfrak{p}\in\operatorname{Spec}(R):\ \mathfrak{p}\supseteq I,\ \mathcal{F}^{E_{S_{\mathfrak{p}}}}\text{ is not finitely generated as a ring over }\widehat{S_{\mathfrak{p}}}\}
\]
is a closed set in the Zariski topology; more precisely, it is equal to $V(J),$ where
\[
J=\bigcap_{F\in\operatorname{IGL}(\Delta)} \mathfrak{p}_F,
\]
$\mathfrak{p}_F$ is the prime ideal generated by the variables whose indices are not in $F,$ and
\[
\operatorname{IGL}(\Delta):=\{F\subset\Delta:\ (K^{[2]}:_R K)\neq K^{[2]}+(\lcm_K),\ K:=(I:_R \mathbf{x}_F)\}.
\]
In particular, $W\cap\hbox{}^*\operatorname{Spec} (R)=\{\mathfrak{p}_F:\ F\in\operatorname{IGL}(\Delta)\}.$
\end{teo}

\begin{proof}
Let $\mathfrak{p}$ be a prime ideal of $R$ with $\mathfrak{p}\supset I$ and such that $\mathcal{A}^{S_{\mathfrak{p}}}$ is a not finitely generated ring over $S_{\mathfrak{p}}.$ Thanks to Lemma \ref{check only with graded prime}, we know that $\mathfrak{p}\supset\mathfrak{p}^*,$ and $\mathcal{A}^{S_{\mathfrak{p}^*}}$ is also a not finitely generated ring over $S_{\mathfrak{p}^*}.$ Now, write $\mathfrak{p}^*=\mathfrak{p}_F$ for some $F\subset [n].$ Since $\mathfrak{p}_F\supset I,$ by \cite[Theorem 1.7]{MillerSturmfels2005} there is a minimal prime ideal $\mathfrak{p}_G$ of $I$ such that $\mathfrak{p}_F\supset\mathfrak{p}_G$ for some $G\subset\Delta.$ This implies that $F\subset G$ and therefore, since $G\subset\Delta$ and $\Delta$ is a simplicial complex, $F\subset\Delta.$ This shows that $\mathfrak{p}\supset J,$ and therefore $W\subset V(J).$

Now, let $\mathfrak{p}\in V(J),$ so $\mathfrak{p}\supset\mathfrak{p}_F$ for some $F\in\operatorname{IGL}(\Delta).$ By construction, $\mathfrak{p}_F\in W;$ moreover, since $W$ is closed under specialization by Proposition \ref{the principal locus is stable under generalization} one has that $\mathfrak{p}\in W,$ just what we finally wanted to prove.
\end{proof}

Theorem \ref{main result: calculation of the infinitely generated locus} leads to a very naive algorithm for computing the infinitely generated locus of the Frobenius algebra attached to a Stanley--Reisner ring; in this method, our input is a simplicial complex $\Delta$ as before and our output will be the ideal $J$ described along the statement of Theorem \ref{main result: calculation of the infinitely generated locus}. This method has been implemented in Macaulay2 \cite{M2}.

\begin{enumerate}[(i)]

\item\textbf{Step 1}: Initialize $K$ as $I,$ $\mathfrak{p}$ as the ideal $(x_1,\ldots,x_n),$ and $L$ as the empty list.

\item\textbf{Step 2}: If $(K^{[2]}:K)\neq K^{[2]}+\langle\lcm_K\rangle,$ then add $\mathfrak{p}$ to the list $L.$ Otherwise, return the empty list.

\item\textbf{Step 3}: For each non--empty face $F\subset\Delta,$ assign to $K$ the value $(I:\mathbf{x}_{F}).$ If $(K^{[2]}:K)\neq K^{[2]}+\langle\lcm_K\rangle,$ then add $\mathfrak{p}=\mathfrak{p}_F$ to the list $L,$ where $\mathfrak{p}$ denotes the ideal generated by the variables that do not belong to $F.$

\item\textbf{Step 4}: Output the intersection of the elements of $L.$

\end{enumerate}

\begin{disc}\label{our algorithm just looking at the simplicial complex}
Thanks to \cite[Theorem 4]{AlvarezYanagawa2014}, we know that the Frobenius algebra attached to the injective hull of a complete Stanley--Reisner ring is finitely generated if and only if $\Delta$ has no free faces. In this way, we can rewrite the above algorithm in terms of the simplicial complex $\Delta$ as follows.

\begin{enumerate}[(i)]

\item If $\Delta$ has no free faces, then stop and output the empty list.

\item For each non--empty face $F\subset\Delta,$ let $\link(F)$ be its corresponding link. If $\link (F)$ has at least one free face, then add $\mathfrak{p}=\mathfrak{p}_F$ to the list $L.$

\end{enumerate}

\end{disc}
We conclude this section by exhibiting a family of ideals with a quite simple infinitely generated locus; before so, we want to review the following notion \cite[Definition 2.5]{BoocherSeinerBEH}.

\begin{df}\label{nearly complete intersection ideals: definition}
Let $\mathbb{K}$ be any field, let $R=\mathbb{K}[x_1,\ldots,x_n]$ be the polynomial ring in $n$ variables over $\mathbb{K},$ let $I\subseteq R$ be a squarefree monomial ideal, and let $\supp (I)$ be the support of $I;$ that is,
\[
\supp (I):=\{i\in [n]:\ x_i\text{ divides at least one minimal monomial generator of } I\}.
\]
We say that $I$ is a \textbf{nearly complete intersection} (hereafter, NCI for short) if it is generated in degree at least two, is not a complete intersection, and for each $i\in\supp(I),$ the monomial localization $I(\mathfrak{p}_{([n]\setminus\supp (I))\cup\{i\}})$ is a complete intersection.
\end{df}
The interested reader in this family of ideals may also like to consult \cite{MillerStoneNCI} where a classification of these ideals in the degree two case is given. Our reason for considering this class of ideals is given in the following:

\begin{prop}\label{infinitely generated locus in the NCI case}
Let $\K$ be a field of prime characteristic $p,$ let $R=\K [x_1,\ldots, x_n],$ $S=R/I,$ and $I$ is a NCI ideal. Then, the set
\[
W=\{\mathfrak{p}\in\operatorname{Spec}(R):\ \mathfrak{p}\supseteq I,\ \mathcal{A}^{S_{\mathfrak{p}}}\text{ is not finitely generated as a ring over }S_{\mathfrak{p}}\}
\]
is either empty or is given by the closed set $V(J),$ where $J:=(x_i:\ i\in\supp (I)).$
\end{prop}

\begin{proof}
If $W=\emptyset,$ then we are done, so we assume hereafter $W\neq\emptyset.$ First of all, notice that
\[
J=(x_i:\ i\in\supp (I))=\mathfrak{p}_{[n]\setminus\supp (I)}
\]
is a prime ideal. Secondly, since $W\neq\emptyset,$ in particular one has that $J\in W.$ Now, we assume, to reach a contradiction, that there is a face ideal $\mathfrak{p}\in W$ such that $\mathfrak{p}\subsetneq J.$ This implies that there is $i\in\supp (I)$ such that
\[
\mathfrak{p}\subseteq (x_j:\ j\in\supp(I)\setminus\{i\})=\mathfrak{p}_{([n]\setminus\supp (I))\cup\{i\}}.
\]
However, since $W$ is closed under specialization and $\mathfrak{p}\in W,$ we have that $\mathfrak{p}_{([n]\setminus\supp (I))\cup\{i\}}\in W,$ which is equivalent to say, thanks to Lemma \ref{links and localized Frobenius algebras}, that
\[
(I'^{[2]}: I')\neq I'^{[2]}+(\lcm_{I'}),
\]
where $I'=I(\mathfrak{p}_{([n]\setminus\supp (I))\cup\{i\}}).$ But this is a contradiction because, since $I$ is NCI by assumption, $I'$ is a complete intersection and therefore in this case we know that
\[
(I'^{[2]}: I')= I'^{[2]}+(\lcm_{I'}).
\]
Summing up, we have checked that, if $W\neq\emptyset,$ then
\[
W=V(\mathfrak{p}_{[n]\setminus\supp (I)})=V((x_i:\ i\in\supp (I))),
\]
as claimed.
\end{proof}




\section{Examples}\label{Examples}
The goal of this section is to show several examples of how to use our method in Macaulay2; the first one is directly borrowed from \cite[Example 5]{BoixZarzuela2017}.

\begin{ex}
Let $\Delta$ be the simplicial complex given by facets $\{1,2,3\},\ \{1,2,6\}$ and $\{3,4,5\}.$ In the following calculation, the set of vertices $\{1,2,3,4,5,6\}$ is identified with the set of variables $\{x,y,z,w,a,b\}.$

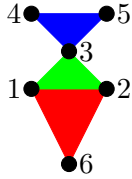
\smallskip

\begin{center}
\begin{tikzpicture}[scale=1]
\draw[fill, red] (0.5,0)--(0,1)--(1,1)--(0.5,0);
\draw[fill, green] (0,1)--(0.5,1.5)--(1,1)--(0,1);
\draw[fill,blue] (0.5,1.5)--(0,2)--(1,2)--(0.5,1.5);
\draw[fill] (0.5,0) circle [radius=0.1];
\draw[fill] (0,1) circle [radius=0.1];
\draw[fill] (1,1) circle [radius=0.1];
\draw[fill] (0.5,1.5) circle [radius=0.1];
\draw[fill] (0,2) circle [radius=0.1];
\draw[fill] (1,2) circle [radius=0.1];
\node[right] at (0.5,0) {6};
\node[left] at (0,1) {1};
\node[right] at (1,1) {2};
\node[right] at (0.5,1.5) {3};
\node[left] at (0,2) {4};
\node[right] at (1,2) {5};
\end{tikzpicture}
\end{center}

We use our method to determine the non--finitely generated locus of the corresponding Stanley--Reisner ideal.
\begin{verbatim}

loadPackage "SimplicialComplexes";

load "non_finitely_generated_locus_algorithm.m2";

R=QQ[x,y,z,w,a,b];

A= simplicialComplex {x*y*z,x*y*b,z*w*a};

I=monomialIdeal(A);

I

monomialIdeal (x*w, y*w, x*a, y*a, z*b, w*b, a*b)

nonfglocus(I)

monomialIdeal (x, y, w, a, b)
\end{verbatim}
Notice that $\link (3)$ is a subsimplicial complex of $\Delta$ with facets $\{1,2\}$ and $\{4,5\},$ and its corresponding squarefree monomial ideal has an infinitely generated Frobenius algebra \cite[Example 1]{BoixZarzuela2017}. Here, we have depicted this link:

\smallskip 

\begin{center}
\begin{tikzpicture}[scale=1]
\draw (0,0)--(1,0);
\draw (0,1)--(1,1);
\draw[fill] (0,0) circle [radius=0.1];
\draw[fill] (1,0) circle [radius=0.1];
\draw[fill] (0,1) circle [radius=0.1];
\draw[fill] (1,1) circle [radius=0.1];
\node[left] at (0,0) {1};
\node[right] at (1,0) {2};
\node[left] at (0,1) {4};
\node[right] at (1,1) {5};
\end{tikzpicture}
\end{center}
\end{ex}

The below example is exactly \cite[Example 2.10]{Gallegoetal2021}.

\begin{ex}
Let $\Delta$ be the simplicial complex given by facets $\{1,2,5\},\ \{1,3,5\}$ and $\{1,2,4\}.$

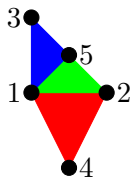
\smallskip

\begin{center}
\begin{tikzpicture}[scale=1]
\draw[fill, red] (0.5,0)--(0,1)--(1,1)--(0.5,0);
\draw[fill, green] (0,1)--(0.5,1.5)--(1,1)--(0,1);
\draw[fill,blue]
(0,1)--(0,2)--(0.5,1.5)--(0,1);
\draw[fill] (0.5,0) circle [radius=0.1];
\draw[fill] (0,1) circle [radius=0.1];
\draw[fill] (1,1) circle [radius=0.1];
\draw[fill] (0.5,1.5) circle [radius=0.1];
\draw[fill] (0,2) circle [radius=0.1];
\node[right] at (0.5,0) {4};
\node[left] at (0,1) {1};
\node[right] at (1,1) {2};
\node[right] at (0.5,1.5) {5};
\node[left] at (0,2) {3};
\end{tikzpicture}
\end{center}

We use our method to determine the non--finitely generated locus of the corresponding Stanley--Reisner ideal.
\begin{verbatim}

R=QQ[x_1,x_2,x_3,x_4,x_5];

A= simplicialComplex {x_1*x_2*x_5,x_1*x_3*x_5,x_1*x_2*x_4};

I=monomialIdeal(A);

I

monomialIdeal (x x , x x , x x )
                2 3   3 4   4 5

nonfglocus(I)

monomialIdeal (x , x , x , x )
                2   3   4   5
\end{verbatim}
\end{ex}
We can also explore how far are from Theorem \ref{main result: calculation of the infinitely generated locus} some of the components calculated in \cite{Gallegoetal2021} using our method in the following way. The first example is the one studied in \cite[Proposition 2.6]{Gallegoetal2021}.

\begin{ex}
Let $\Delta$ be the simplicial complex given by facets $\{2\}$ and $\{1,3\}.$ As it is explained there, the upper bound given by \cite[Theorem 2.4]{Gallegoetal2021} yields the ideal $(x_2).$ Our algorithm shows, in this particular example, that the whole non--finitely generated locus is given by a smaller ideal.
\begin{verbatim}

R=QQ[x_1,x_2,x_3];

I=monomialIdeal(x_1*x_2,x_2*x_3);

{closedcomponent1(I),nonfglocus(I)}

{monomialIdeal x , monomialIdeal (x , x , x )}
                2                  1   2   3
\end{verbatim}
\end{ex}
The next example is the one studied in \cite[Example 2.9]{Gallegoetal2021}. This example shows, in particular, that the upper bound given by \cite[Theorem 2.7]{Gallegoetal2021} is, in general, not always equal to the defining ideal of the non--finitely generated locus.

\begin{ex}
Let $\Delta$ be the simplicial complex given by facets $\{1,2,4\},\ \{1,3\}$ and $\{2,3\}.$ We proceed as before.
\begin{verbatim}

R=QQ[x_1,x_2,x_3,x_4];

I=monomialIdeal(x_1*x_2*x_3,x_3*x_4);

{closedcomponent1(I),nonfglocus(I)}

{monomialIdeal x , monomialIdeal (x x , x , x )}
                3                  1 2   3   4
\end{verbatim}
\end{ex}

\section*{Acknowledgments}
The authors would like to specially thank the referee for his/her useful comments and remarks concerning the contents of this paper.

Alberto F. Boix and S. Zarzuela are partially supported by Spanish Ministerio de Econom\'ia y Competitividad grant PID2019-104844GB-I00.  D. A.
J. G\'omez-Ram\'irez thanks the Instituci\'on Universitaria Pascual Bravo,
Visi\'on Real Cognitiva S.A.S; and Juan Carlos Díaz and Liliana Lopera for all his kindness and support.

\bibliographystyle{alpha}
\bibliography{AFBoixReferences}

\end{document}